\documentclass[a4paper,12pt, legno]{article}
\usepackage{amsmath}
\usepackage{amsfonts}
\usepackage{amsthm}
\usepackage{mathrsfs}
\usepackage[T1]{fontenc}
\usepackage{amsmath}
\usepackage{amsfonts}
\usepackage{amsthm}
\usepackage{mathrsfs}
\usepackage{amssymb}
\usepackage{pst-node}
\usepackage{tikz-cd}

\newtheorem{cor}{Corollary}

\newtheorem{rem}[]{Remark}

\newcommand{\Rset}{\mathbb{R}}
\newcommand{\Cset}{\mathbb{C}}

\begin{document}
December 15, 2021.

\begin{center}
\textbf{ON BACKGROUND DRIVING DISTRIBUTION FUNCTIONS (BDDF) FOR SOME SELFDECOMPOSABLE VARIABLES.}
\end{center}
\begin{center}
Zbigniew J. JUREK\footnote{Institute of Mathematics, University od Wroclaw, Pl. Grunwaldzki 2/4,

\noindent  50-384 Wroclaw, Poland;
\ \ zjjurek@math.uni.wroc.pl; \ \  www.math.uni.wroc.pl/$^\sim$zjjurek}
\end{center}
\medskip
[for  final version see \emph{Mathematica Applicanda}  vol. 49(2),2021, pp.85-104.]

\medskip
\medskip
\emph{Abstract.} Many classical variables (statistics) are selfdecomposable. They admit the random integral representations via L\'evy processes. In this note are given formulas for their background driving distribution functions (BDDF). This may be used for a simulation of those variables. Among the examples discussed are: gamma variables, hyperbolic characteristic functions, Student t-distributions, stochastic area under planar Brownian motions, inverse Gaussian variable, logistic distributions, non-central chi-square, Bessel densities and Fisher z-distributions. Found representations
might be of use in statistical applications.



hyperbolic distribution; t-Student distribution; logistic distribution;

non-central chi-square distribution; Wiener process; stochastic area

variable;  Bessel functions;  Fisher z-distribution.

\medskip
\medskip
In the probability or statistic course Gaussian distribution and the central limit theorem (CLT) occupy the main stage. Often as an extension
stable distributions are introduced and then the class of infinitely divisible distributions (ID). But there is a subclass of selfdecomposable distributions also known as the class $L$. Surprisingly many distributions used in statistics and other areas of applications  are in $L$. Class $L$ distributions have a random integral representation which maybe useful in simulations. For many classical distributions in statistics  we find explicite formulas for their background driving distribution (BDDF). To have a better insight in those distributions (how they are concentrated) there are also  computed some their numerical values using the platform WolframAlpha.com. In Jurek and Kepczynski (2021)will be given the graphs for  the BDDF computed in Section \textbf{B)} below.

\newpage
\textbf{A). Selfdecomposability:  basic formulas  and relationships.}

\medskip
If a sequence  $ X_k, k=1,2,...$  of random variables is \emph{strongly mixing}, in particular, if it is stochastically independent, $a_k, k=1,2,...$ is a sequence of positive numbers such that
the triangular array $a_nX_k; 1\le k \le n, n\ge 1$ is infinitesimal and $b_k, k=1,2 , ... $ is a sequence of real numbers such that
\begin{equation}
a_n(X_1+X_2+...+X_n) + b_n \Rightarrow X, \ \mbox{in the weak topology, as}\ \ n\to \infty, \ \ \ \
\end{equation}
then $X$ is called \emph{selfdecomposable} or \emph{L\'evy class L variable}.

Those variables are characterized as follows:
\begin{equation}
X\in L \ \mbox{iff} \ \  \forall(0<c<1)\exists \ (X_c\, \emph{idependent of} \, X) \   X \stackrel{d}{=}cX+X_c .  \ \ \ \ \ \  \ \ \ \
\end{equation}
See  Bradley and Jurek (2014) for the case of strong mixing sequences. For the independent sequences and some historical origins of the problem cf.  Gnedenko and Kolmogorov (1954), Section 29-30, or Loeve (1963), Section 23 or Feller (1966), Chapter XVII.8 . (The above property is used to justify the terminology \emph{selfdecomposability, selfdecomposable}.)

\medskip
\noindent Furthermore, from Jurek and Vervaat (1983) we have the following \emph{random integral  representation} for selfdecomposable variables:
\emph{
\begin{multline*}
X \mbox{ is a selfdecomposable} \ \mbox{iff} \\ \mbox{there exists unique (in distribution)  L\'evy process $(Y_X(t), t\ge 0)$ such that} \\
 X=\int_0^\infty e^{-t}dY_X(t), \ \mbox{and} \  \mathbb{E}[\log(1+|Y_X(1)|)]<\infty, \ \ \ \ \  \ \ \ \  (\star)
\end{multline*} }
If the identity $(\star)$ holds then to the variable $Y_X(1)$ we refer to as the  BDRV (\emph{background driving random variable}) of $X$. So, $Y_X(1)$ is infinitely divisible with finite log-moment; in short,   $Y_X(1)\in ID_{\log}$.

Similarly, if
$\phi_X$ and $\psi_X $ are the characteristic  functions of $X$ and $Y_X(1)$, respectively, then for $t\neq 0$
$$\psi_X(t)\equiv\psi_{Y_X(1)}(t)=  \exp [t (\log \phi_{X}(t))^\prime] =\exp [t \frac{\phi_{X}^{\prime}(t)}{\phi_X(t)}] \in ID_{\log}, \ \    \ \ \                     \ (\star  \star). $$
Therefore, if the identity $(\star \star)$ holds then  to the function $\psi_X$ we will refer as to the BDCF (\emph{ the background driving characteristic function}) of $X$; cf. Jurek (2001), Proposition 3.

We say that $X \in L$ has \emph{the factorization property} if $X$ convoluted  with its BDCF is in $L$ again, that is,
\begin{equation}
 \ \ \ \ \ \ \ \ \ \ \ L^f:= \{X \in L: \phi_X(t)\cdot \psi_X(t)\in L\}\subsetneq L \subset ID; \ \ \ \ \ \ \ \  \ \ \ \
\end{equation}
cf. Jurek (1985), Theorem 4.5 , Corollary 4.6  and  Iksanov, Jurek and Schreiber (2004), Theorem 1, Proposition 3; more on the class $L^f$  is in  Czyżewska and Jurek (2011).

Finally, if $G_X(a):=P(Y_X(1)\le a)$ is the BDDF (\emph{background driving distribution function}) of $X$ then it is given by the following formula:
\begin{equation*}
G_X(a)= \frac{1}{2}-\frac{1}{\pi}\int_0^\infty \Im(\exp\big[ -ita+t \frac{\phi_{X}^\prime(t)}{\phi_X(t)}\big])\frac{dt}{t}, \ \ a\in C_{G_X}, \ \ \ \  \ \ \ (\star \star \star)
\end{equation*}
where $C_{G_X}$ is the set of continuity  points the probability distribution function $G_{X}\in ID_{\log}$ (infinitely divisible with finite log-moment); cf. Jurek (2019), Proposition 1 and references therein.

\medskip
\noindent If the characteristic function $\phi_X(t)$ is a real one then
$$
G_X(a)= \frac{1}{2}+\frac{1}{\pi}\int_0^\infty \exp\big[ t\frac{\phi_{X}^{\prime}(t)}{\phi_{X}(t)}\big]\frac{\sin(ta)}{t}dt, \ a\in C_{G_X}, \ \ \  \ \  \ \ (\star \star \star \star)
$$
is the BDDF of $X$ and $G_X(-a)=1-G_X(a), a \in C_{G_X}$ .

\medskip
\underline{\textbf{Symbolically we may summarize  all the above as:}}
\begin{equation}
X \in L \ \ \mbox{iff} \ \ X=\int_0^\infty e^{-s}dY_X(s), \  \ P(Y_X(1)\le a)=G_X(a) \in ID_{log}, \ \  \ \
\end{equation}
where $a\in C_{G_X}$ and  $(Y_X(s), s\ge 0),$ is the unique (in distribution) L\'evy process for $X$.  \ \ \ \ \ \

\medskip
\medskip
Finally, the class L forms a closed (in the weak topology) convolutions semigroup of $ID$ (all infinitely divisible variables). If $X\in L$ then $aX+b\in L$, for all reals $a$ and $b$. All $X\in L$ are absolutely continuous  with respect to Lebesque measure, i.e., $X\in L$ have probability densities. All non-generate class $L\subsetneq ID$ distributions have infinite  L\'evy (spectral) measures in their L\'evy-Kintchine representations .

\medskip
\noindent \textbf{\underline{ The following selfdecomposable distribution are discussed below:}}

\medskip
(1) gamma;  (2) chi-square; (3) log-gamma;  (4) inverse gamma;

(5) hyperbolic-cosine; (6) hyperbolic-sine;   (7) hyperbolic-tangent;

(8) via zeros of Bessel function;

(9) t-Student;

(10) stochastic area under planar Brownian motion;

(11) exponential integral and inverse Gaussian;

(12) logistic distribution;

(13)  non-central chi-square and Bessel $h_{\nu}$ distribution;

(14) Fisher z-distribution.

\medskip
\textbf{Warning.} All  numerical values below computed via WollframAlpha  are approximate not exact one,  although the equality sign is used.

\newpage
\textbf{B). Examples of  BDDF for some selfdecomposable distributions.}

\medskip
\textbf{1).}  $X:=\gamma_{\alpha,\lambda}$ - \textbf{the gamma random variable.}

\medskip
Let us recall that $f_{\gamma_{\alpha,\lambda}}(x)=\frac{\lambda^\alpha}{\Gamma(\alpha)}x^{\alpha-1}e^{-\lambda x}, \ \ x\in (0,\infty),$ is the probability density function ($\alpha>0$ is the \emph{shape parameter} and $\lambda>0$ is the \emph{rate parameter})  and
$\phi_{\gamma_{\alpha, \lambda}}(t)=(1- it/\lambda)^{-\alpha}$ is its characteristic function which is selfdecomposable by Jurek (1997), p. 97. From $(\star \star)$ we get the BDCF
$\psi_{\gamma_{\alpha , \lambda}}(t)= \exp\alpha\,[\frac{1}{1-it/\lambda}-1]$ (compound Poisson distribution).
Finally we have (three formulae for) the  BDDF:
\begin{multline}
G_{\gamma_{\alpha,\lambda}}(a)=P\big( \sum_{k=1}^{N_{\alpha}}\mathcal{E}_k (\lambda) \le a \big)\\ = \frac{1}{2}+\frac{1}{\pi}\int_0^\infty \exp(-\frac{\alpha\,(t/\lambda)^2}{1+ (t/\lambda)^2})\, \sin \big(t a- \frac{\alpha\,(t/\lambda)}{1+(t/\lambda)^2}\big)\,\frac{1}{t}dt \\
 = e^{-\alpha}+ e^{-\alpha} \int_0^{2\sqrt{\alpha  \lambda a}}I_1(w)e^{-w^2/4\alpha}dw, \  \ \ a\in C_{G_X}\cap (0,\infty);  \qquad \qquad
\end{multline}
where $N_\alpha$ is Poisson variable with the parameter $\alpha$ independent of i.i.d. exponential $\mathcal{E}_k$ and  finally  $I_1(x)$ is the (modified)  Bessel function; cf. Jurek (2019), Lemma 1 and the identity $(\star \star \star)$. Thus from (4) we may write that

$$\gamma_{\alpha, \lambda}=\int_0^\infty e^{-t}dY_{\gamma_{\alpha, \lambda}}(t); \  \ \ P(Y_{\gamma_{\alpha,\lambda}}(1)\le a)=G_{\gamma_{\alpha,\lambda}}(a) \in ID_{\log} .$$

\medskip
\underline{\textbf{Here are some numerical values:}}

$ G_{\gamma_{2,1}}(0.001)=0.135606 $; \    $G_{\gamma_{2,1}}(0.01)=0.138042$; \  $G_{\gamma_{2,1}}(1)=0.394297 $; \ $G_{\gamma_{2,1}}(2)=0.6035$;   \ $ G_{\gamma_{2,1}}(3)=0.753011$; \  $G_{\gamma_{2,1}}(4)= 0.8519$; \  $ G_{\gamma_{2,1}}(6)=0.95123$.

\medskip
\medskip
\medskip
\textbf{2). $X:=\chi^2(n)\equiv \gamma_{n/2,1/2}$  -  chi-square distributions.}

\medskip
Recall that $X\equiv\chi^2(n):= \sum_{m=1}^n Z_m^2, \ \mbox{where} \  \ Z_m \ \mbox{are i. i. d. $N(0,1)$}$ is called \emph{the chi-square distribution with the $n$ degrees of freedom}. Then
$\phi_{\chi^2(n)}(t)= \frac{1}{(1-2it)^{n/2}}\in L$ and its BDCF is  $\psi_{\chi^2(n)}(t)=\exp n/2[\frac{1}{1-i 2t}-1] $. Hence by $(\star \star \star) $ we have
\begin{equation}
G_{\chi^2(n)}(a)= \frac{1}{2}+\frac{1}{\pi}\int_0^\infty \exp(-  \frac{2nt^2}{1+4t^2})\, \sin (at - \frac{nt}{1+4t^2}) \frac{dt}{t}, \ \mbox{for} \ a>0.
\end{equation}
So, there exists  a L\'evy process $(Y_{\chi^2(n)}(t), t\ge 0)$ such that $Y_{\chi^2(n)}(1)\stackrel{d}{=}G_{\chi^2(n)}$ and
$$\chi^2(n)=\int_0^\infty e^{-t}dY_{\chi^2(n)}(t), \ \ P(Y_{\chi^2(n)}(1)\le a)= G_{\chi^2(n)}(a)\in ID_{\log}.  $$

 \medskip
\textbf{Note.} In particular, for  $\chi^2(2)=N(0,1)^2+\tilde{N}(0,1)^2$ (a sum of squares of two independent standard normal rv) we have
$$
G_{\chi^2(2)}(a)= \frac{1}{2} +\frac{1}{\pi}\int_0^\infty \exp[-4t^2/(1+4t^2)]\sin (ta-2t/(1+4t^2))dt/t.
$$

\medskip
\underline{\textbf{Here are some numerical values:}}

 $ G_{\chi^2(2)}(1)=0.53013, \    G_{\chi^2(2)}(3)=0.7477, \  G_{\chi^2(2)}(5)=0.8686,$  \   \

 $ G_{\chi^2(2)}(7)= 0.9332, \ G_{\chi^2(2)}(10)=0.9766, \ G_{\chi^2(2)}(15)= 0.9962. $

\medskip
\medskip
\textbf{3). $ X:= \log \gamma_{\alpha,\lambda}$ is the log-gamma variable.}

\medskip
Note that $f_{\log \gamma_{\alpha,1}}(x)= \frac{1}{\Gamma(\alpha)}e^{\alpha x - e^x}, x \in \Rset$ is the probability density of $\log \gamma_{\alpha,1}$. Hence $\log \gamma_{\alpha, \lambda}= -\log \lambda +\log \gamma_{\alpha,1}$. Furthermore,
from Jurek (1997),Example c), p. 98, we  have that $\log \gamma_{\alpha,\lambda} \in L$ (is selfecomposable) and its has the characteristic function is
$\phi_{\log \gamma_{\alpha,\lambda}}(t)=\lambda^{-it}\frac{\Gamma(\alpha+it)}{\Gamma(\alpha)}$.

Using $(\star \star)$, log-gamma variables has the BDCF  of the form
$$
\psi_{\log \gamma_{\alpha,\lambda}}(t))= \exp( - it \log\lambda+it \psi(\alpha+it),
$$
where on the right-hand side the  $\psi(z):= \frac{d}{dz}\log\Gamma(z)$ denotes \emph{ the digamma function}. Finally, by $(\star \star \star)$,
\begin{equation}
G_{\log \gamma_{\alpha, \lambda}}(a)=\frac{1}{2}-\frac{1}{\pi}\int_0^\infty \Im(e^{-ita -it\log\lambda+it\psi(\alpha+it)})\frac{dt}{t};
\end{equation}
is the BDDF of log-gamma variable. Thus we have
$$
\log\gamma_{\alpha,\lambda}=\int_0^\infty e^{-s}dY_{\log\gamma_{\alpha,\lambda}}(s), \ \ P(Y_{\log\gamma_{\alpha,\lambda}}(1)\le a)=G_{\log \gamma_{\alpha, \lambda}}(a)\in ID_{\log} .
$$

\medskip
\medskip
\underline{\textbf{Here are some numerical values:}}

$G_{\log \gamma_{2,1}}(-2)= 0.03, \ \ G_{\log \gamma_{2, 1}}(-1)= 0.109, \ \  G_{\log \gamma_{2,1}}(0)= 0.3099;$

$ G_{\log \gamma_{2, 1}}(1)=0.6635;  \ \ G_{\log \gamma_{2, 1}}(2)=0.9503. $

\medskip
\medskip
\medskip
For the references below let us recall that functions
\begin{equation}
I_{\nu}(z)=\sum_{k=0}^\infty \frac{(\frac{1}{2} z)^{2k+\nu}}{k!\Gamma(\nu+k+1)}, \ \ \ K_n(z)=\frac{\pi}{2}\frac{I_{-n}(z)- I_n(z)}{\sin(\pi n)}, \ \nu \in \Rset, \  \ z \in \Cset,
\end{equation}
are called \emph{the modified Bessel functions of the first and second kind}, respectively. These are solutions to the Bessel second order differential equation; cf. Gradshteyn and Ryzhik (1994), Sections 8.40, 8.43 and the formula \textbf{8.445}.

\medskip
\medskip
\textbf{4). $ X:= 1/\gamma_{\alpha,\lambda}$ is the inverse - gamma variable.}

\medskip
From Jurek(2001), Proposition 1 with Example 1 (c) we have that random variable $1/\gamma_{\alpha,\lambda}$ is selfdecomposable and
\begin{multline*}
f_{1/\gamma_{\alpha,\lambda}}(x)=\frac{\lambda^\alpha}{\Gamma(\alpha)} (\frac{1}{x})^{\alpha+1} e^{-\lambda/x}, x>0, \ \ \mbox{ is its pdf and}, \\     \phi_{1/\gamma_{\alpha,\lambda}}(t)=\frac{2}{\Gamma(\alpha)} (-i\lambda t)^{\alpha/2}K_{\alpha}(2\sqrt{-i \lambda t}), t\in \Rset,  \ \ \mbox{ is the characteristic function;} \ \ \
\end{multline*}
($K_{\alpha}(z)$ denotes the modified Bessel function). Thus the BDDF for the inverse-gamma variable is
\begin{multline}
G_{1/\gamma_{\alpha,\lambda}}(a)=\frac{1}{2} -\frac{1}{\pi}\int_0^\infty \Im\big(exp[-ita -\sqrt{-i\lambda t}\, \frac{K_{\alpha-1}(2\sqrt{-i\lambda t})}{K_{\alpha}(2\sqrt{-i\lambda t})}]\big)\frac{1}{t}dt.
\end{multline}

\medskip
\noindent This so, because by Wolframalpha,
$$
(\log K_{\alpha}(2\sqrt{-i\lambda t}))^\prime= \frac{i\lambda}{2\sqrt{-i \lambda t}} \frac{K_{\alpha-1}(2\sqrt{-i\lambda t})+ K_{\alpha+1}(2 \sqrt{-i\lambda t})}{K_{\alpha}(2\sqrt{-i\lambda t})},
$$
and hence we get
$$
t(\log \phi_{1/\gamma_{\alpha,\lambda}}(t))\prime= \frac{\alpha}{2}-\frac{\sqrt{-i\lambda t}}{2}\, \frac{K_{\alpha-1}(2\sqrt{-i\lambda t})+ K_{\alpha+1}(2 \sqrt{-i\lambda t})}{K_{\alpha}(2\sqrt{-i\lambda t})}.
$$
\noindent Again by Wolframalpha we have  the identity
 $$  \frac{z}{2}\big[\frac{K_{\alpha-1}(z)+K_{\alpha+1}(z)}{K_{\alpha}(z)}]= \alpha+ z \frac{K_{\alpha-1}(z)}{K_\alpha(z)}.$$
Applying it for  $z=2 \sqrt{-i \lambda t}$ we get, by $(\star \star)$, the BDCF formula $\psi_{1/\gamma_{\alpha, \lambda}}(t)$ and then we infer the
formula (9).

Consequently, we have
$$
1/\gamma_{\alpha,\lambda}=\int_0^\infty e^{-s}dY_{1/\gamma_{\alpha,\lambda}}(s),  \ \ P(Y_{1/\gamma_{\alpha,\lambda}}(1) \le a)=G_{1/\gamma_{\alpha,\lambda}}(a)\in ID_{\log}.
$$

\medskip
\underline{\textbf{Numerical illustrations for $\alpha=\lambda=2$.}}

$ G_{1/\gamma_{2,2}}(0)= 3.807\times 10^{-11}; $
$ G_{1/\gamma_{2,2}}(1/10)= 0.00318;  $  $ G_{1/\gamma_{2,2}}(2/10)= 0.0501;  $ $ G_{1/\gamma_{2,2}}(0.5)= 0.292043; $
$ G_{1/\gamma_{2,2}}(1)=0.550257; $ $ G_{1/\gamma_{2,2}}(2)=0.7645; $ $ G_{1/\gamma_{2,2}}(3)= 0.851994;  $ $ G_{1/\gamma_{2,2}}(5)=0.924258;  $      $ G_{1/\gamma_{2,2}}(6)=0.941699; $ $ G_{1/\gamma_{2,2}}(10)=0.973368; $ $ G_{1/\gamma_{2,2}}(20)= 0.992454;  $

\medskip
\medskip
\textbf{EXAMPLE 1.} Let $W_t, t\ge 0$ be a Wiener process in $\Rset^d, d \ge 3$ (the escape process), starting from zero .Then
$$
L_r:=\sup\{t\ge 0; ||W_t||\le r\}, r\ge 0, \ \ L_r\stackrel{d}{=} \frac{1}{2\gamma_{\frac{d-2}{2}, r^2}},
$$
and $L_r, r\ge 0$ is a process with independent increments (but not homogenous for $d>3$) and continuous in probability,cf. Getoor (1979).

\medskip
\medskip
\textbf{5). X: = $\hat{C}$ - hyperbolic-cosine variable};

\medskip
A random varible $\hat{C}:=\frac{2}{\pi}\sum_{k=1}^\infty \frac{1}{2k-1}\eta_k$ , (almost surely converging series),
where $\eta_k$ are i.i.d. standard Laplace variables has a the characteristic function
 $\phi_{\hat{C}}(t)= \frac{1}{\cosh t}\in L$, cf. Jurek (1996). Its BDCF  is $\psi_{\hat{C}}(t)=\exp(- t \tanh t)$ and by
 $(\star \star \star \star )$ its BDDF is given by
\begin{equation}
G_{\hat{C}}(a)= \frac{1}{2}+\frac{1}{\pi}\int_0^\infty \exp(-x \tanh x)\,\frac{\sin (ax)}{x}dx.
\end{equation}
By the identity (4) we have the identification:
$$
\hat{C}=\int_0^\infty e^{-s}dY_{\hat{C}}(s), \ \ P(Y_{\hat{C}}(1)\le a)= G_{\hat{C}}(a)\in ID_{\log}.
$$
\medskip
\medskip
\textbf{6). X: = $\hat{S}$ - hyperbolic-sine variable};

\medskip
Note that $\hat{S}:=\frac{1}{\pi}\sum_{k=1}^\infty \frac{1}{k}\eta_k$, where $\eta_k$, as above in \textbf{5)}, are i.i.d. Laplace r.v. and its characteristic  function  is
$\phi_{\hat{S}}(t):= \frac{t}{\sinh t}$;  (symmetric variable), cf. Jurek (1996). By $(\star \star )$ its BDCF  is equal to $\psi_{\hat{S}}(t)=\exp(1- t\coth t)$ that finally, by $(\star \star \star \star)$, gives BDDF for $\hat{S}$ as follows
\begin{equation}
G_{\hat{S}}(a)= \frac{1}{2}+\frac{1}{\pi}\int_0^\infty \exp(1 -x \coth x)\,\frac{\sin (ax)}{x}dx.
\end{equation}
In view of (4) we have
$$
\hat{S}=\int_0^\infty e^{-t}dY_{\hat{S}}(t), \ \ P(Y_{\hat{S}}(1)\le a)= G_{\hat{S}}(a)\in ID_{\log}.
$$
\medskip
\emph{
\begin{rem}
a) Since $\Gamma(1+it)\Gamma(1-it)=\frac{\pi t}{\sinh(\pi t)}$ by GR 8.332 (3), from section \textbf{3).} above, we infer
$
\pi \hat{S}\stackrel{d}{=}\log\mathcal{E}_1 - \log \mathcal{E}_2=\log \frac{\mathcal{E}_1}{\mathcal{E}_2},$
where $\mathcal{E}_i, i=1,2 ,$  are i.i.d. standard exponential variables.
\newline
b) Note that the (normalized) hyperbolic-sine
$\frac{2}{\pi^2}\frac{x}{\sinh x}, x\in \Rset,$
is also  a probability density function and its characteristic function is
$$\int_{-\infty}^\infty e^{itx}\frac{2}{\pi^2}\frac{x}{\sinh x}dx= 4\frac{e^{\pi t}}{(1+e^{\pi t})^2}=\frac{1}{\cosh^2(\pi t/2)} \in L,$$
equivalently \  $\int_{0}^\infty e^{itx}\frac{x}{\sinh(\pi x/2)}dx=\frac{1}{\cosh^2(t)}$;
\noindent by \textbf{GR} page 1185(23); (note that there  the Fourier transform  is  $1/\sqrt{2\pi}$ times the characteristic function !!)
\end{rem}
}

\medskip
\medskip
\medskip
\textbf{7). X := $\hat{T}$ - hyperbolic-tangent};

Recall that afunction  $\phi_{\hat{T}}(t)= \frac{\tanh t}{t}$ is called \emph{ the hyperbolic-tangent characteristic function}. We have that
$\hat{T}\in L$ (is selfdecomposable) by Jurek (1996). By $(\star \star)$ it has BDCF $\psi_{\hat{T}}(t)=\exp[\frac{2t}{\sinh(2t)}-1]\in ID_{\log}$. It is a compound Poisson distribution $\sum_{k=1}^{N_1}\xi_k)$, where $\xi_k\stackrel{d}{=} \hat{C} \ast \hat{S} $ where $\hat{C}$ and $\hat{S}$ are independent hyperbolic-cosine and sine variables as  $\frac{2t}{\sinh(2t)}=\frac{1}{\cosh t}\frac{t}{\sinh t}$).

\medskip
From $(\star \star \star \star)$ we have that
\begin{equation}
G_{\hat{T}}(a)=\frac{1}{2}+\frac{1}{\pi}\int_0^\infty\exp[\frac{2t}{\sinh(2t)}-1]\frac{\sin(ta)}{t}dt,
\end{equation}
is the BDDF of hyperbolic-tangent random variable $\hat{T}$. Furthermore, by (4) we have
$$
\hat{T}=\int_0^\infty e^{-t}dY_{\hat{T}}(t), \ \ P(Y_{\hat{T}}(1))\le a)= G_{\hat{T}}(a).
$$

\medskip
\underline{\textbf{Here are some numerical values:}}

$G_T(0.4)=0.7653; G_T(1)=0.8645; G_T(2)=0.9528; G_T(3)=0.9846.$

\medskip
\medskip
\medskip
\textbf{Example 2.} For a Brownian Motion $(B(t),t\ge 0)$, let $\tau_1$ be the exit time of BM from the interval $[-1,1]$, i.e., $\tau_1:=\inf\{t:|B(t)|=1\}$.

If  $g_{\tau_1}:=\sup\{t<\tau_1: B(t)=0\}$ is  the last zero of $(B(t), t\ge 0)$ before exiting $[-1,1]$ then $B(g_{\tau_1})\stackrel{d}{=} N\sqrt{g_{\tau_{1}}}\in L$, where $N$ is independent normal variable. Moreover, $\phi_{B(g_{\tau_1})}(t)=\tanh(t)/t$;
cf. Yor (1997), Sect.18.6,  p. 133; also Jurek (2001), Example 1(b).

From the above formula  and the identity $(\star)$, we conclude that there is a L\'evy process $(Y_{B(g_{\tau_1})}(t), t\ge 0)$ such that $Y_{B(g_{\tau_1})}(1)\stackrel{d}{=}G_{\hat{T}}$ and
$$
B(g_{\tau_1})=\int_0^\infty e^{-s}dY_{B(g(\tau_{1}))}(s),
$$
where $P(Y_{B(g(\tau_{1}))}(1) \le a)= \frac{1}{2}+\frac{1}{\pi}\int_0^\infty\exp[\frac{2t}{\sinh(2t)}-1]\frac{\sin(ta)}{t}dt$.

\medskip
\medskip
\emph{
\begin{rem}For the hyperbolic characteristic functions we have the obvious equality: $\phi_{\hat{C}}(t)=\phi_{\hat{S}}(t)\cdot \phi_{\hat{T}}(t)$. Hence we infer that on the level of their corresponding BDCF we have $\psi_{\hat{C}}(t)=\psi_{\hat{S}}(t)\cdot \psi_{\hat{T}}(t)$.
Finally, from (10),(11) and (12) we get $(1-t \coth t)+(\frac{2t}{\sinh (2t)}-1)= -t \tanh t$, that is, the formula for $G_{\hat{C}}$ can be obtained for those of $G_{\hat{S}}$ and $G_{\hat{T}}$.
\end{rem}
}

\medskip
\medskip
\medskip
\textbf{8). Class L distributions via zeros of Bessel function.}

\medskip
From Jurek (2003), Theorem 1, the converging series
$$
X_{\nu}:= \sum_{k=1}^\infty z^{-1}_{\nu,k}\eta_k \in L, \ \mbox{where $\eta_k$ are i.i.d. Laplace rv and}$$
and $ z_{\nu,k}\ \mbox{are zeros of the function} \ z^{-\nu}J_{\nu}(z)$ is well defined.
The variable $X_{\nu}$ has  the characteristic  function $\phi_{X_{\nu}}(t)= (2^{\nu}\Gamma (\nu+1))^{-1} \, \frac{t^{\nu}}{I_{\nu}(t)}\in \Rset$. Consequently, by $(\star \star)$, its BDCF is
\begin{multline*}
\psi_{X_\nu}(t)=\exp[t(\log\phi_{X_\nu}(t))^\prime]= \exp[t(\nu/t - (\log I_\nu(t))^\prime)]\\
=\exp[\nu - t(I_{\nu-1}(t)/I_\nu(t)-\nu/t)=\exp[2\nu -  tI_{\nu-1}(t)/I_\nu(t)].
\end{multline*}
Hence, by $(\star \star \star \star)$,
\begin{equation}
G_{X_\nu}(a)=\frac{1}{2}+\frac{1}{\pi}\int_0^\infty\exp[2\nu -  tI_{\nu-1}(t)/I_\nu(t)]\frac{\sin(at)}{t}dt
\end{equation}
is BDDF for $X_\nu$ and there is a L\'evy process $(Y_{X_\nu}(t), t\ge 0)$ such that
$$
X_\nu=\int_0^\infty e^{-t}dY_{X_\nu}(t), \ \ P(Y_{X_{\nu}}(1)\le a) = G_{X_\nu}(a)\in ID_{\log}.
$$

\medskip
\medskip
\underline{\textbf{For an illustration here  are some numerical values:}}

$G_{X_2}(0.001)=0.500757, \ \ G_{X_2}(0.01)= 0.5075, \   G_{X_2}(0.1)=0.57,$  \

$G_{X_2}(0.5)=0.82,$ \ $ G_{X_2}(1)=0.95, \  G_{X_2}(2)=0.998,\ G_{X_2}(3)=0.999973.$

\medskip
\medskip
\textbf{9). $X:=T^{\nu}$ - Student t- distributions, $\nu>0$.}

\medskip
A random variable $ X:=T^{\nu}, \nu >0$ with the probability density function
$$
f_{\nu}(x)=\frac{\Gamma(\nu+1/2)}{\sqrt{2\pi\nu}\,\Gamma(\nu)}\,(1+\frac{x^2}{2\nu})^{-\nu-1/2}, \ \ x \in\Rset,
$$
is called \emph{the Student t-distribution with $2\nu$ degrees of freedom}. It  has the characteristic function
$$
\phi_{T^{\nu}}(t)= \frac{2^{1-\nu}}{\Gamma(\nu)}\,(\sqrt{2\nu}\,|t|)^\nu K_{\nu}(\sqrt{2\nu}\,|t|)\in L,
$$
cf.Jurek (2001), Example 2, p. 247. The BDCF for $T^{\nu}$ is equal to:
\begin{equation*}
\psi_{T^{\nu}}(t) = \exp\Big[ -|t|\sqrt{2\nu}\,\frac{K_{\nu-1}(\sqrt{2\nu}|t|)}{K_{\nu}(\sqrt{2\nu}|t|)} \Big], \ \ t\neq 0,
\end{equation*}
where $K_{\nu}$ is the modified Bessel function of the second kind; cf. the end of Section \textbf{3).} Thus from $(\star \star \star \star)$
\begin{equation}
G_{Y_{T^{\nu}}(1)}(a) = \frac{1}{2}+\frac{1}{\pi}\int_0^\infty \exp \Big[ -|t|\sqrt{2\nu}\,\frac{K_{\nu-1}(\sqrt{2\nu}|t|)}{K_{\nu}(\sqrt{2\nu}|t|)} \Big]\frac{\sin(ta)}{t}dt.
\end{equation}
Finally appealing to (4) we get
$$
T^{\nu}=\int_0^\infty e^{-t}dY_{T^{\nu}}(t),  \ \ \ P(Y_{T^{\nu}}(1) \le a)= G_{T^{\nu}}(a) \in ID_{\log}.
$$

\medskip
\medskip
\underline{\textbf{Here are some numerical values} for $\nu=2$}:

$ G_{Y_{T^{2}}(1)}(0.02)=0.50558, \ \ G_{Y_{T^{2}}(1)}(0.5)=0.6253, \ \  G_{Y_{T^{2}}(1)}(1)=0.7458,$
$ G_{Y_{T^{2}}(1)}(2)=0.8888; $\ \ $ G_{Y_{T^{2}}(1)}(3)=0.9497, \ \ G_{Y_{T^{2}}(1)}(4)=0.9756,$ \ \

$G_{Y_{T^{2}}(1)}(10)=0.9988 $

\medskip
\medskip
\emph{
\begin{rem}
Since $K_{\nu}(z)/K_{-\nu}(z)=1$, therefore for Student t-distribution $T^{1/2}$, with one degree of freedom, we recover its BDDF as
$$G_{Y_{T^{1/2}}(1)}(a)= \frac{1}{2}+\frac{1}{\pi}\int_0^\infty e^{-t}\frac{\sin ta}{t}dt = \frac{1}{2}+\frac{1}{\pi}\arctan a . $$
And that is the 1-stable distribution. This is not a surprise as stable distributions are fixed points of the random integral mapping $(\star)$, cf. Jurek and Vervaat (1983), Theorem 5.1.
\end{rem}
}

\medskip
\medskip
\textbf{ 10). The stochastic area under the planar Brownian motion.}

\medskip
 \textbf{a).  \ $\mathcal{A}_u:=\int_0^uZ_sd\tilde{Z}_s-\tilde{Z}_sdZ_s, u>0. $}

 \medskip
 Let $(Z_t, \tilde{Z}_t)$ be the planar Brownian Motion (BM) and we define $\mathcal{A}_u:=\int_0^uZ_sd\tilde{Z}_s-\tilde{Z}_sdZ_s, u>0.$  Then the stochastic area $\mathcal{A}_1$, under the graph given that $(Z_1, \tilde{Z}_1)=(1,1)$ has the conditional characteristic function
$$
\phi_{\mathcal{A}_1}(t):=\mathbb{E}[e^{i t \mathcal{A}_1}|(Z_1,\tilde{Z}_1)=(1,1)]=\frac{t}{\sinh t}\exp[ -(t \coth t -1)]
$$
by P. L\'evy (1950) or Yor (1992), p. 19.

Since $t/\sinh t \in L $ has the  factorization property (cf. (3) in section \textbf{A})  and $\exp[-t(\coth t-1)]$ is its BDCF, therefore  by  Iksanov, Jurek, Schreiber (2004), Proposition 1 and 3), the product $ \frac{t}{\sinh t}\exp[ -(t \coth t -1)]$, is a selfdecomosable characteristic function . By $(\star \star)$ its BDCF  is

$
\exp t( \log ( \frac{t}{\sinh t}\exp[ -(t \coth t -1)] ) )\prime = \exp t( \log t-\log\sinh t- t\coth t +1 )^\prime =$
$
\exp t( 1/t -\coth t -\coth t + t/\sinh^2 t)=\exp(1-2 t\coth t +t^2/\sinh^2t).  \  \
$

Consequently, by  the formula$(\star \star \star \star)$, in Section \textbf{A)} above, we get
\begin{multline}
G_{\mathcal{A}_1}(a):= P[\mathbb{A}_1\le a| (Z_1, \tilde{Z}_1)=(1,1)]\\=\frac{1}{2}+ \frac{1}{\pi}\int_0^\infty \exp[1- 2 t\coth t+ t^2/\sinh^2 t]
\frac{\sin (ta)}{t}dt.
\end{multline}
Finally, by (4) we have that there is a L\'evy process such that
$$
\mathcal{A}_1 = \int_0^\infty e^{-t}dY_{\mathcal{A}_1}(t), \ \ \ \ P(Y_{\mathcal{A}_1}(1)\le a)=G_{\mathcal{A}_1}(a)\in ID_{\log}.
$$

\medskip
\medskip
\underline{\textbf{Here are some  numerical values:}}

$G_{\mathcal{A}_1}(0.5)=0.649892; \   G_{\mathcal{A}_1}(1)=0.775697, \ \  G_{\mathcal{A}_1}(1.2)=0.8163$ \ \

$G_{\mathcal{A}_1}(1.5)= 0.86674; \
G_{\mathcal{A}_1}(2)=0.92558; \ G_{\mathcal{A}_1}(3)= 0.9799$

\medskip
\medskip
\medskip
\textbf{b). \  $\mathbb{A}_u^p:=\int_0^u V^p_sd\tilde{V}^p_s-\tilde{V}^p_sdV^p_s, u>0; \ \ p> -1/2.$}

\medskip
\noindent This represents the generalized L\'evy stochastic area for the process:

$\textbf{V}^p_t=(V^p_t,\tilde{V}^p_t) :=t^{-p}\int_0^t s^p d\textbf{B}_s$, where $(\textbf{B}_s, s\ge0)$ is the planar Brownian  motion. Note that for $p=0$ we get the case \textbf{a)}, above.

Assuming that
$\textbf{V}^p_1=(V^p_1,\tilde{V}^p_1)=(1,1)$ then for the conditional stochastic area $\mathbb{A}_1^p$ we have
\begin{equation*}
\phi_{\mathbb{A}_1^p}(t):= \mathbb{E}[e^{it \mathbb{A}^p_1}|\mathbf{V}^p_1=(1,1)]=(2^\nu \Gamma(\nu+1))^{-1} \frac{|t|^\nu}{I_{\nu}(|t|)}\exp\big[- |t| \frac{I_{\nu+1}(|t|)}{I_{\nu}(|t|)}\big]
\end{equation*}
where $I_\nu(z)$ (cf. the end of Section \textbf{3).}) is the modified Bessel function and $\nu:=p+1/2>0$; cf. Biane and Yor (1987) or Yor (1989) or Duplantier (1989); also Jurek (2003).

\medskip
Since $(2^\nu \Gamma(\nu+1))^{-1} \frac{|t|^\nu}{I_{\nu}(|t|)}\in L_f\subset L$ (see Iksanov, Jurek, Schreiber (2004), Propositions 1 and 3), then  using Theorem 1 in Jurek (2003) we have that $\phi_{\mathbb{A}_1^p}(t)\in L.$ Hence

$$
\log \phi_{\mathbb{A}^p_1}(t)= - \log(2^\nu \Gamma(\nu+1)) + \nu \log t -\log I_{\nu}(t) - t \ \frac{I_{\nu+1}(t)}{I_{\nu}(t)}.
$$
From Wolframalpha.com we have the identities:
 \begin{multline*}
  (i)\ (\log I_{\nu}(t))^{\prime} = \frac{I_{\nu-1}(t)}{I_{\nu}(t)} -\frac{\nu}{t};  \  \ \ \  (ii)\ \  2\nu - t \frac{I_{\nu-1}(t)+I_{\nu+1}(t)}{I_{\nu}(t)}=-2t \frac{I_{\nu+1}(t)}{I_\nu(t)}  \\ (iii) \  \ (\frac{I_{\nu+1}(t)}{I_\nu(t)})^\prime=
 \frac{1}{2}\big( 1+ I_{\nu+2}(t)/I_{\nu}(t)-(I_{\nu+1}(t)/I_{\nu}(t))^2-I_{\nu-1}(t)I_{\nu+1}(t)/I_{\nu}^2(t) \big)  \ \ \ \ \ \           \ \ \ \ \
 \end{multline*}
Hence we get the derivative:
\begin{multline*}
(\log \phi_{\mathbb{A}^p_1}(t))^\prime =\nu/t - (\frac{I_{\nu-1}(t)}{I_{\nu}(t)} -\frac{\nu}{t})- (I_{\nu+1}/I_{\nu} +t (I_{\nu+1}/I_{\nu})^\prime))\\
=2\nu/t-I_{\nu-1}/I_{\nu}-I_{\nu+1}/I_{\nu}-t(I_{\nu+1}/I_{\nu})^\prime \\ =2\nu/t-I_{\nu-1}/I_{\nu}-I_{\nu+1}/I_{\nu}-t/2\big( 1+ I_{\nu+2}/I_{\nu}-(I_{\nu+1}/I_{\nu})^2-I_{\nu-1}I_{\nu+1}/I_{\nu}^2 \big),
\end{multline*}
and consequently by $(\star \star)$ we get the BDCF
\begin{multline*}
\psi_{\mathbb{A}^p_1}(t)=\exp[t (\log\phi_X(t))^\prime]\\ = \exp[ 2\nu - t \frac{I_{\nu-1}(t)+I_{\nu+1}(t)}{I_{\nu}(t)} -\frac{t^2}{2}(1+\frac{I_{\nu+2}(t)}{I_{\nu}(t)}- (\frac{I_{\nu+1}(t)}{I_{\nu}(t)})^2 - \frac{I_{\nu-1}(t)I_{\nu+1}(t)}{I_{\nu}(t)^2}]\\ =
\exp[ -2 t \frac{I_{\nu+1}(t)}{I_{\nu}(t)} -\frac{t^2}{2}(1+\frac{I_{\nu+2}(t)}{I_{\nu}(t)}- (\frac{I_{\nu+1}(t)}{I_{\nu}(t)})^2 - \frac{I_{\nu-1}(t)I_{\nu+1}(t)}{I_{\nu}(t)^2}],
\end{multline*}
where we used  the identity (ii). Hence we have  the formula for BDDF as follows:
\begin{multline}
G_{{\mathbb{A}^p_1}}(a)= \frac{1}{2}+ \frac{1}{\pi} \int_0^\infty \Big( \exp[- 2t \frac{I_{\nu+1}(t)}{I_{\nu}(t)}]\\ \exp[-\frac{t^2}{2}\big(1+\frac{I_{\nu+2}(t)}{I_{\nu}(t)}- (\frac{I_{\nu+1}(t)}{I_{\nu}(t)})^2 - \frac{I_{\nu-1}(t)I_{\nu+1}(t)}{I_{\nu}(t)^2}\big)]\Big)\frac{\sin(ta)}{t}dt.
\end{multline}
From the identification (4), from section \textbf{A)} we have
$$
\mathbb{A}^p_1 = \int_0^\infty e^{-t}dY_{\mathbb{A}^p_1}(t) , \ \ \ P(Y_{\mathbb{A}^p_1}(1)\le a)=G_{{\mathbb{A}^p_1}}(a)\in ID_{\log}.
$$

\medskip
\medskip
\underline{\textbf{Some numerical values for} $\nu=2 \equiv p=3/2$}:

$G_{\mathbb{A}^{3/2}_1}(0.1)=0.5420; \ \  G_{\mathbb{A}^{3/2}_1}(0.3)=0.6239; \ \   , G_{\mathbb{A}^{3/2}_1}(0.5)=0.700;$ \ \

$G_{\mathbb{A}^{3/2}_1}(1)=0.84422; \ \    G_{\mathbb{A}^{3/2}_1}(2)= 0.97553; \ \  G_{\mathbb{A}^{3/2}_1}(3)=0.997414$.

\medskip
\emph{
\begin{rem}
By some (tedious) calculation one may check directly  that  we have equal two  BDDF $G_{\mathcal{A}_1}= G_{\mathbb{A}_1^0}$ , i.e, by putting $\nu=1/2$ in (16) we get equality (14).
Similarly, we have  the equality on the level of their corresponding characteristic functions $\phi_{\mathcal{A}_1}(t)=\phi_{\mathbb{A}_1^0}(t), t \in \Rset.$
\end{rem}
}

\medskip
\medskip
\textbf{11). Some integral functionals of Wiener Process.}

\medskip
\textbf{a). Exponential integrals.}

\medskip
For the process $(W_t+ \alpha t, t\ge 0)$,  (a parameter  $\alpha >0$)  and the variable $X:=\int_0^\infty e^{- W_s - \alpha t}ds$, Urbanik (1992) in Example 3.3 (on p.309), proved
$$
X\stackrel{d}{=}1/\gamma_{2\alpha,2}, \ \mbox{which has pdf} \ \ \frac{4^\alpha}{\Gamma(2 \alpha)}x^{-2\alpha-1}e^{-2/x}, x>0
$$
Thus
$$
\int_0^\infty e^{- W_t - \alpha t}dt\stackrel{d}{=}\int_0^\infty e^{-s}dY_{1/\gamma_{2\alpha,2}}(s), \ \ Y_{1/\gamma_{2\alpha,2}}(1)\stackrel{d}{=}G_{1/\gamma_{2\alpha,2}}.
$$
The  BDDF for the inverse gamma was given in \textbf{4).} above.

\medskip
\medskip
\emph{\begin{rem}
Note that if we define
$\tau_s:=\inf (t>0: W_t+\alpha t>s)$  then  $(\tau_s,s\ge 0)$ is the inverse Gaussian (L\'evy) process .
\end{rem} }

\medskip
\medskip
\textbf{ b). Inverse Gaussian  variable $IG(\lambda,\mu)$} .

For $\lambda>0,\mu>0,$ random variables with probability density functions $f_{IG(\lambda,\mu)}$ and characteristic functions $\phi_{IG(\lambda,\mu)}$ given as
$$
f_{IG(\lambda,\mu)}(x)=\sqrt{\frac{\lambda}{2\pi x^3}} \exp[-\frac{\lambda(x-\mu)^2}{2\mu^2 x}]; \ \ \ \phi_{IG(\lambda,\mu)}(t)=\exp\Big[\frac{\lambda}{\mu}\big(1-\sqrt{1-\frac{2\mu^2 it}{\lambda}}\big)\Big ],
$$
are called \emph{the inverse Gaussian distributions}. By Halgreen (1979) we have that they are selfdecomposable.
Since
$ (\log \phi_{IG(\lambda,\mu)}(t))^\prime = i \mu (1-\frac{2\mu^2it}{\lambda})^{-1/2}$,
thus its BDCF is
$\psi_{IG(\lambda,\mu)}(t)= \exp[i\mu t (1-2\mu^2\lambda^{-1}it)^{-1/2}]$.
 Consequently,
\begin{equation}
G_{IG(\lambda,\mu)}(a)=\frac{1}{2}- \frac{1}{\pi}\int_0^\infty \Im(\exp[-ita + i\mu t \sqrt{1-\frac{2\mu^2it}{\lambda}}\,\,])\frac{dt}{t},
\end{equation}
is the BDDF for $IG(\lambda,\mu)$. So from the identification (4) we get
$$
IG(\lambda,\mu)=\int_0^\infty e^{-t}dY_{IG(\lambda,\mu)}(t),  \ \ \ P(Y_{IG(\lambda,\mu)}(1)\le a)= G_{IG(\lambda,\mu)}(a)\in ID_{\log}.
$$
\medskip
\underline{\textbf{Here are some numerical values for} $\lambda=\mu=1$}:

$G_{IG(1,1)}(-5)=0.00; $  \  $G_{IG(1,1)}(-3)=0.04 ;$ \    $G_{IG(1,1)}(-2)=0.23; $

$G_{IG(1,1)}(-1)=0.55; $ \ $G_{IG(1,1)}(-0.1)=0.77; $ \  $G_{IG(1,1)}(0)=0.79;$ \

$G_{IG(1,1)}(0.1)=0.81 ;$  \  $G_{IG(1,1)}(0.5)= 0,87 ;$ \   $G_{IG(1,1)}(1)= 0.91$;

 $G_{IG(1,1)}(2)= 0.96$;  \ \ $G_{IG(1,1)}(3)=0.98 ;$  $G_{IG(1,1)}(5)=0.99 $

\medskip
\underline{\textbf{Similarly we have for} $\mu=1$ and $\lambda=2$}:

$G_{IG(2,1)}(-3)=0.007$; \  $ G_{IG(2,1)}(-2)=0.14$; \   $G_{IG(2,1)}(-1)=0.54$; \

$G_{IG(2,1)}(-0.5)= 0.72$; \  \ $G_{IG(2,1)}(0)=085$;  \ $G_{IG(2,1)}(0.5)= 0.926$;

\  $G_{IG(2,1)}(1)=0.9638$, \ $G_{IG(2,1)}(2)=0.9914$ .

\medskip
\medskip
\textbf{c). Quadratic functional of BM.}

For the Brownian process  $(W_t, t\ge 0)$ starting from zero, the drift parameter $b$,  the staring starting point $a$  and an independent of $W$ a standard normal variable $N$, the variable
$$
Q(a,b)\equiv Q(a,b,W):= N(\int_0^1(W_s + bs+ a)^2ds)^{1/2} \ \mbox{is selfdecomposable.}
$$
To see this, let us take two Brownian motions $W$ and $\tilde{W}$ and two standard normal variables $N$ and $\tilde{N}$  that are all  independent. Then for $0<c<1$ we have $N\stackrel{d}{=}cN+\sqrt{1-c^2}\tilde{N}$ and hence
\begin{multline}
Q(a,b)=N(\int_0^1(W_s + bs+ a)^2ds)^{1/2} \stackrel{d}{=}(cN+ \sqrt{1-c^2}\tilde{N})(\int_0^1(W_s + bs+ a)^2ds)^{1/2} \\ \stackrel{d}{=}cQ+ \sqrt{1-c^2}\tilde{N}(\int_0^1(\tilde{W}_s + bs+ a)^2ds)^{1/2}\stackrel{d}{=}cQ+\sqrt{1-c^2}\tilde{Q},
\end{multline}
with $Q$ and $\tilde{Q}$ independent variables.This, with the characterization (2),  proves the selfdecomposability property of $Q(a,b)$.

\medskip
From by Wenocur (1986) or from Yor(1992), p.19, formula (2.9) we have the characteristic function
\begin{multline*}
\phi_{Q(a,b)}(t)\\ = \mathbb{E}[\exp( it N(\int_0^1(W_s + bs+ a)^2ds)^{1/2})]= \mathbb{E}[\exp (-\frac{t^2}{2}\int_0^1(W_s+bs+a )^2ds)]\\
= (\cosh t)^{-1/2}\exp[-\frac{b^2}{2}(1- \frac{\tanh t}{t})-ba(1-\frac{1}{\cosh t})-\frac{a^2}{2} t \tanh t)]\\
= (\cosh t)^{-1/2}\exp[-\frac{a^2}{2} t \tanh t)] \,\, \exp[\frac{b^2}{2}(\frac{\tanh t}{t}-1)]\,\,\exp[ba(\frac{1}{\cosh t}-1)].
\end{multline*}

\begin{rem}
\emph{In the  formula above for $\phi_{Q(a,b)}$  we have:
\newline(i) $(\cosh t)^{-1/2} \in L$ has the factorization property ( see (3) in the section \textbf{A)})  and $\exp[-\frac{1}{2} t \tanh t]$ is its
BDCF (which in fact, is s-selfdecomposable distribution; cf. Iksanov,Jurek and Schreiber (2004), p.1367).
\newline(ii) $\exp[\frac{b^2}{2}(\frac{\tanh t}{t}-1)]$ and $\exp[ba(\frac{1}{\cosh t}-1)]$ are compound Poisson of the form $\sum_{k=1}^{N_{\alpha}}\xi_k$. In the first case, $\alpha=b^2/2$ and $\xi_k$ are i.i.d. with hyperbolic tangent variables and in the second case, $\alpha =ba>0$ and $\xi_k$ are i.i.d. hyperbolic cosine variables.}
\end{rem}
We have that
\begin{multline*}
(\log\phi_Q(t))^\prime = (-\frac{1}{2} \log \cosh t- \frac{a^2}{2} t \tanh t+ \frac{b^2}{2}(\frac{\tanh t}{t} -1) +ab (\frac{1}{\cosh }-1))^{\prime}  \\=  -\frac{1}{2} \tanh t- \frac{a^2}{2}(\tanh t +t\frac{1}{\cosh^2t})+\frac{b^2}{2}(\frac{1}{\cosh^2t}t^{-1}-t^{-2}\tanh t)\\- ab(\cosh t)^{-2}\sinh t \\=
-\frac{1}{2}\tanh t(1+a^2+b^2t^{-2})+\frac{1}{2}\frac{1}{\cosh^2t}(-a^2t+ b^2 t^{-1}-2ab \sinh t).
\end{multline*}
Hence, the BDCF is
$$
\psi_{Q(a,b)}(t)=\exp[-\frac{1}{2}\tanh t((1+a^2)t+b^2t^{-1})+\frac{1}{2}\frac{1}{\cosh^2t}(-a^2t^2+ b^2 -2ab t\sinh t) ]
$$
and finally by $(\star\star \star \star)$ we get its BDDF
\begin{multline}
G_{Q(a,b)}(x)=\frac{1}{2} +\frac{1}{\pi}\int_0^\infty\exp[-\frac{1}{2}\tanh t((1+a^2)t+b^2t^{-1})\\+\frac{1}{2}\frac{1}{\cosh^2t}(-a^2t^2+ b^2 -2ab t\sinh t)]\frac{\sin(t x)}{t} dt.
\end{multline}
In summary , by (4) in Section \textbf{A)},  we have that
$$
Q(a,b)=\int_0^\infty e^{-s}dY_{Q(a,b)}(s), \ \ \ P(Y_{Q(a,b)}(1)\le u )=G_{Q(a,b)}(u)\in ID_{\log}
$$

\medskip
\underline{\textbf{Here are some numerical values for} $a=1$ and $b=2$}:

\noindent $Q(1,2)(0.01)=0.501664,  \ Q(1,2)(0.1)= 0.516603,  \ Q(1,2)( 1)=0.648221$, \

\noindent $Q(1,2)(2)=0.763609, \ Q(1,2)(3)= 0.849722, \ \ Q(1,2)(4)= 0.908518$,\

$Q(1,2)(5)=0.966382.$

\medskip
\medskip
\textbf{12). Logistic distribution } \emph{l(a,b)} .

From Ushakov (1999), p. 298 (or W.Feller p. 52) random variable with the probability distribution function
$$
p(x;a,b)=  \frac{\pi}{b\sqrt{3}} \frac{\exp[- \frac{\pi}{\sqrt{3}} (\frac{x-a}{b})]}{\big(1+\exp[-\frac{\pi}{\sqrt{3}}(\frac{x-a}{b})]\big)^2}; \ -\infty<x<\infty \  \ \
( a\in \Rset, \ b>0);
$$
which gives the mean-value $a$, the variance $b^2$,
and  the characteristic  function
\begin{equation*}
\phi_{l(a,b)}(t)=e^{ita}\Gamma(1- ict)\Gamma(1+ict), \ \ \mbox{where} \   c:= b\sqrt{3}/\pi,
\end{equation*}
is called \emph{the logistic distribution}.  It is selfdecomposable by \textbf{3).} and the fact that $L$ is a convolution semigroup; cf. last paragraph in the  section \textbf{A)}.
Since
$
\Gamma(1- ict)\Gamma(1+ict) = \pi ct/\sinh(\pi ct),\ t\in\Rset,
$
by Gradsteyn-Ryzhik \textbf{8332}.3 we get
$$\phi_{l(a,b)}(t)=e^{ita}\, \pi c t/\sinh(\pi ct).$$
Thus by $(\star \star)$, the  BDCF for the logistic $\emph{l}(a,b)$ distribution is
$$
\psi_{\emph{l(a,b)}}(t)=\exp[ita + 1 - \pi c t \coth (\pi ct)].
$$
Finally, the BDDF for \emph{l(a,b)} is:
\begin{equation}
G_{\emph{l(a,b)}}(u)= \frac{1}{2}-\frac{1}{\pi}\int_0^\infty \Im(\exp(-itu +itu + 1 -\pi ct \coth (\pi ct)))\frac{1}{t}dt, \ u \in \Rset.
\end{equation}
Comp. with the hyperbolic-sine function above in \textbf{6).}

\medskip
By  the identification formula (4) we have
$$
\emph{l(a,b)}=\int_0^\infty e^{-t}dY_{\emph{l(a,b)}}(t), \ \ \ \  P( Y_{\emph{l(a,b)}}(1)\le u )= G_{\emph{l(a,b)}}(u) \in ID_{\log}.
$$
\medskip
\medskip
\textbf{Here are some values for} $\emph{l}(0,1)$:

$G_{\emph{l}(0,1)}(0.5)=0.58, \  G_{\emph{l}(0,1)}(1)=0.62, \  G_{\emph{l}(0,1)}(2)=0.8, \ \ G_{\emph{l}(0,1)}(3)=0.89,  \ \ G_{\emph{l}(0,1)}l(5)=0.97$.

\medskip
\medskip
\textbf{13).The Bessel function $I_{\nu}(x)$ appearing in probability density functions.}

\medskip
\medskip
\textbf{(a). Non-central chi-square distribution.}

\medskip
\emph{The non-central chi-square distribution} $\chi_{c}^2(k)$, with \emph{k-degrees of freedom} (k>0) and\emph{the non-centrality exponent $c$} ($c>0$), has the probability density
$$ \frac{1}{2}e^{-(x+c)/2} \Big(\sqrt{\frac{x}{\lambda}}\Big)^{k/2-1} I_{k/2-1}(\sqrt{c x}), \ \  x>0; $$
This is a Poisson mixture of the gamma densities $\gamma_{a_k,1}$  for specified  shape parameters $a_k$ and scales equal 1;
see Feller (1966), p. 58, formula (7.3).)

By WolframAlpha.com
$
\mathbb{E}[\chi_{c}^2(k)]= c+k), \ \ Var[\chi_{c}^2(k)]= 2(2c+k)$, and
the  characteristic  function  $\phi_{\chi_{c}^2(k)}$ is of the form
$$
\phi_{\chi_{c}^2(k)}(t)=\frac{\exp(\frac{itc}{1-2it})}{(1-2it)^{k/2}}= \exp\Big[\frac{c}{2}\big(\frac{1}{1-2it} -1\big)\Big]\, \, \frac{1}{(1-2it)^{k/2}}\in ID.
$$
It corresponds to the sum of two independent ID variable: the first one is the compound Poissona with the Levy-Khintchine formula
$$
\exp[\int_0^\infty(e^{itx}-1)c \frac{1}{2} e^{-x/2}dx]
$$
(not in L  as its Levy (spectral) measure is finite !) and the second one is the chi-square with $k/2$ degrees of freedom and it is in L,(see \textbf{2)} above)  with its Levy-Khintchine-formula
$$
\frac{1}{(1-2it)^{k/2}}= \exp[\int_0^\infty(e^{itx}-1)\frac{k}{2} \frac{e^{-x/2}}{x}dx].
$$
All in all, the non-central $\chi_{c}^2(k)$ -distribution is infinitely divisible and its L\'evy (spectral) measure

$$M(dx):= c \frac{1}{2} e^{-x/2}dx + \frac{k}{2} \frac{e^{-x/2}}{x}dx = [\frac{1}{2}e^{-x/2}(c+k/x)]dx $$

has a density such that the function $x[\frac{1}{2}e^{-x/2}(c+k/x)]$ is decreasing on the positive half-line. Consequently,
the non-central $\chi_{c}^2(k)$-distribution is selfdecomposable; Jurek and Mason (1993), p. 94, or  Jurek (1997), p. 96.

Since
\begin{multline*}
(\log (\frac{\exp(\frac{it\lambda}{1-2it})}{(1-2it)^{k/2}}))^{\prime}= (\lambda/2)(\frac{1}{1-2it}-1)^{\prime} -k/2(\log (1-2it))^{\prime}\\ =\lambda/2\frac{2i}{(1-2it)^2} -k/2 \frac{-2i}{1-2it}= i (\frac{\lambda}{(1-2it)^2}+\frac{k}{1-2it}),
\end{multline*}
from $(\star \star)$ we have  the BDCF
$\psi_{\chi_{c}^2(k)}(t)= \exp[it(\frac{c}{(1-2it)^2}+\frac{k}{1-2it})]$, and
finally, the BDDF of $\chi_{c}^2(k)$ is
\begin{equation}
G_{\chi_{c}^2(k)}(a)=\frac{1}{2}-\frac{1}{\pi}\int_0^\infty \Im (\exp[-ita + it(\frac{c}{(1-2it)^2}+\frac{k}{1-2it})])\frac{dt}{t}
\end{equation}

\medskip
\medskip
\underline{\textbf{Here are some numerical values} when $(c=1,k=2)$ } AND the above integral is \textbf{from t=0 to t=10}!!

$   G_{\chi_{1}^2(2)}(1)= 0.4729;    \ \ \ G_{\chi_{1}^2(2)}(2)= 0.55157; \ \ \ G_{\chi_{1}^2(2)}(4)= 0.709$ ;  \ \ \

$ G_{\chi_{1}^2(2)}(8)= 0.88; \ \ \ G_{\chi_{1}^2(2)}(10)=0.93;  \ \ \  G_{\chi_{1}^2(2)}(15)= 0.98  $

\medskip
\medskip
\textbf{(b). Bessel $h_\nu$ - density. }

\medskip
For $\nu>0$,
$h_{\nu}(x):=e^{-x}\frac{\nu I_{\nu}(x)}{x},0<x<\infty,$ is a probability  density function of random variable $\tilde{h}_\nu$.  (Note that for $\nu=r \in \{1,2,...\}$  it is the probability density of the first passage time of symmetric random walk trough r.)

For the characteristic function for $t \in \Rset$ we have
\begin{multline}
\phi_{\tilde{h}_\nu}(t):=\int_0^\infty e^{itx} h_{\nu}(x)dx =[1-it-\sqrt{(1-it)^2-1}]^\nu \\ =\exp \nu\big[it\int_0^\infty \frac{e^{-x}I_0(x)}{1+x^2}dx + \int_0^\infty\big(e^{itx}-1-\frac{itx}{1+x^2}\big)\frac{e^{-x}I_0(x)}{x}dx\big],
\end{multline}
see  Feller (1966), p.414 and 476  or Ushakov (1999), p. 283, for the first  equality. For the second one see a proof below.

From (22), multiplying the density of the L\'evy (spectral) measure by $x$ we get decreasing function $\nu I_0(x)e^{-x}$. Therefore, similarly  as for the non-central chi-square above, we get that Bessel density $h_{\nu}$ ( the random variable $\tilde{h}_\nu$) is selfdecomposable.

Since $(\log(1-it-\sqrt{(1-it)^2-1}))^\prime= i/\sqrt{-t(t+2i)}$ therefore, by $(\star \star)$, the BDCF for the selfdecomposable Bessel density $h_{\nu}$ is
$$
\psi_{{\tilde{h}_1}}(t)=\exp(it/\sqrt{-t(t+2i)})=\exp i \sqrt{\frac{-t}{t+2i}} =\exp i\sqrt{\frac{2}{2-it}-1} \in ID_{\log}
$$
and by $(\star \star \star)$ the BDDF is given as
\begin{equation}
G_{\tilde{h}_\nu}(a)=\frac{1}{2}-\frac{1}{\pi}\int_0^\infty \Im\big(\,\exp[-ita+ i\nu\sqrt{\frac{1}{1-it/2}-1}\,]\,\big)\frac{dt}{t}, a\in C_{h_\nu}.
\end{equation}
Hence by (4), from section \textbf{A)}, we have
$$
\tilde{h}_{\nu}=\int_0^\infty e^{-t}dY_{\tilde{h}_\nu}(t), \ \ P(Y_{\tilde{h}_\nu}(1)\le a) = G_{\tilde{h}_\nu}(a)\in ID_{\log}.
$$

\emph{
\begin{rem}
Note that if $\mathcal{E}(2)$ denotes the exponential  distribution with the parameter 2 then the compound Poisson $\exp(\mathcal{E}(2))$ has the characteristic function $\exp(\frac{1}{1-it/2}-1)$ and its logarithm is under the square root in (23).
\end{rem}}

\medskip
\textbf{Some numerical values  for}  $G_{\tilde{h}_{10}}(\cdot)$:

 $G_{\tilde{h}_{10}}(1)$= 0.0091 ,
  $G_{\tilde{h}_{10}}(5)$= 0.03430,
 $G_{\tilde{h}_{10}}(8)$=0.09192 ,

 $G_{\tilde{h}_{10}}(10)$= 0.1272 ,
 $G_{\tilde{h}_{10}}(200)$= 0.717063,
 $G_{\tilde{h}_{10}}(900)$= 0.868878,

 $G_{\tilde{h}_{10}}(1200)$= 0.882203,
 $G_{\tilde{h}_{10}}(1500)$=0.899366,
 $G_{\tilde{h}_{10}}(2000)$=0.911412 .

\medskip
\medskip
\underline{\emph{A proof of the second equality in (22).}}

\noindent Since $(\phi_{h_{\nu}}(t))^{1/n}= \phi_{h_{\nu/n}}(t)$ we have that
\[
n(\phi_{h_{\nu/n}}(t)-1)=n ((\phi_{h_{\nu}}(t))^{1/n} -1)  \to \log\phi_{h_{\nu}}(t) , \ \mbox{as} \ n\to \infty,
\]
where $h_{\nu/n}(x)=e^{-x}\frac{\nu/n I_{\nu/n}(x)}{x}, \ \ 0<x< \infty$. Thus
\begin{multline*}
n(\phi_{h_{\nu/n}}(t)-1)=n (\int_{0}^\infty (e^{itx}-1)e^{-x}\frac{\nu/n I_{\nu/n}(x)}{x}dx \\ =  it \int_0^\infty x/(1+x^2)e^{-x}\frac{\nu I_{\nu/n}(x)}{x}dx +  \int_0^\infty (e^{itx}-1-\frac{itx}{1+x^2})\,e^{-x}
\frac{\nu I_{\nu/n}(x)}{x}dx \\ \to it(\nu\int_0^\infty \frac{e^{-x}I_{0}(x)}{1+x^2} dx) + \int_0^\infty (e^{itx}-1-\frac{itx}{1+x^2})\,\nu
\frac{e^{-x} I_{0}(x)}{x}\,]dx,
\end{multline*}
which concludes the proof of the second equality in (22).

From (22) we have two ways to the expression for $t(\log\phi_{h_{\nu}}(t))^\prime$ which leads to the equality:
\begin{cor}
For $t>0$ we have
\begin{multline}
\frac{1}{t}\sqrt{\frac{-t}{t+2i}}=\int_0^\infty\frac{e^{-x}I_0(x)}{1+x^2}dx+\int_0^\infty(e^{itx}-\frac{1}{1+x^2})e^{-x}I_0(x)dx \\ =\int_0^\infty e^{itx} e^{-x}I_0(x)dx. \qquad \qquad \qquad \qquad \qquad
\end{multline}
\end{cor}

\medskip
\medskip
\textbf{14). Fisher z-distribution of the dispersion proportion.}

 For positive integers $m_1\ge 1$ and $m_2\ge 1$, degrees of freedom, \emph{  the $z(m_1,m_2)$
 -distribution} (introduced by R. A. Fisher in 1924) has the probability density
 $$
 f_{z(m_1,m_2)}(x):=\frac{2m_1^{m_1/2}m_2^{m_2/2}\Gamma(\frac{m_1+m_2}{2})}{\Gamma(\frac{m_1}{2})\Gamma(\frac{m_2}{2})}\,\frac{e^{m_1x}}{(m_2+m_1 e^{2x})^{(m_1+m_2)/2}}
 $$
with the expected value $\frac{1}{2}(\frac{1}{m_2}-\frac{1}{m_1})$ and the variance $\frac{1}{2}(\frac{1}{m_1}+\frac{1}{m_2})$,
 and the characteristic function
 \begin{equation}
 \phi_{z(m_1,m_2)}(t)= (\frac{m_2}{m_1})^{it/2}\, \frac{\Gamma(\frac{m_1+it}{2})\Gamma(\frac{m_2-it}{2})}{\Gamma(\frac{m_1}{2})\Gamma(\frac{m_2}{2})},
 \end{equation}
 cf. Ushakov (1999), p.309.
 Changing the  parametrization $\alpha_1:= m_1/2, \alpha_2:= m_2/2$ and  then using formulas from \textbf{3)} we infer
 $$
 \phi_{2\,z(2\alpha_1,2\alpha_2)}(t)= (\frac{\alpha_2}{\alpha_1})^{it}\, \frac{\Gamma(\alpha_1+it)\Gamma(\alpha_2-it)}{\Gamma(\alpha_1)\Gamma(\alpha_2)}=
 \phi_{\log \gamma_{\alpha_1,\alpha_1}}(t)\,\ \phi_{- \log \gamma_{\alpha_2,\alpha_2}}(t)\in L,
 $$
 that is, Fisher z-distributions are selfdecomposable. Furthermore,
 their BDCF are $exp[it(\log(\alpha_2/\alpha_1)+\psi(\alpha_1+it)-\psi(\alpha_2-it))]$.

 Finally, for the BDDF we have
 \begin{equation}
 G_{2z(2\alpha_1,2\alpha_2)}(a)=\frac{1}{2}-\frac{1}{\pi}\int_0^\infty\Im\big(e^{-ita + it(\log(\alpha_2/\alpha_1)+\psi(\alpha_1+it)-\psi(\alpha_2-it))}\big)\frac{dt}{t}.
 \end{equation}
As in previous cases, by (4), we have that
$$
2z(2\alpha_1,2 \alpha_2)=\int_0^\infty e^{-t}dY_{2 z(2\alpha_1,2 \alpha_2)}(t), \ P(Y_{2 z(2\alpha_1,2 \alpha_2)}(1)\le a) =  G_{2z(2\alpha_1,2\alpha_2)}(a)\in IG_{\log}.
$$

\medskip
\medskip
 \underline{\textbf{Here are some values for z- distribution with}} $\alpha_1=1$ and $\alpha_2=2$:
 \begin{multline*}
 G_{2z(2,4)}(-5)=0.02727; G_{2z(2,4)}(-3)= 0.1022; G_{2z(2,4)}(-2)= 0.18881; \\  G_{2z(2,4)}(-1)= 0.3296; G_{2z(2,4)}(-0.01)= 0.522759; G_{2z(2,4)}(0)=0.524879; \\  G_{2z(2,4)}(0.01)= 0.5270; G_{2z(2,4)}(0.1)=0.5461; G_{2z(2,4)}(0.5)=0.63107; \\ G_{2z(2,4)}(1)= 073103; G_{2z(2,4)}(2)= 0.8818;  G_{2z(2,4)}(3)=0.9582; \\ G_{2z(2,4)}(4)=0.987396; G_{2z(2,4)}(5)=0.996587.
 \end{multline*}

\medskip
As a by-product of the above we get a relation between the z-distribution and the log-gamma:
\begin{cor}
For two independent gamma random variables $\gamma_{\alpha_1,\alpha_1}$ and $\tilde{\gamma}_{\alpha_2,\alpha_2}$ we have
$$ 2z(2\alpha_1,2\alpha_2)\stackrel{d}{=} \log \gamma_{\alpha_1,\alpha_1} + (-\log \tilde{\gamma}_{\alpha_2,\alpha_2})= \log\frac{\gamma_{\alpha_1,\alpha_1}}{\tilde{\gamma}_{\alpha_2,\alpha_2}}.$$
\end{cor}

\medskip
\medskip
\emph{
\begin{rem}
The classical Fisher - Snedecor distribution $F(d_1,d_2)\stackrel{d}{=}\frac{\chi^2(d_1)/d_1}{\chi^2(d_2)/d_2}$ (a ratio of two independent chi-square distributions with  degrees of freedom $d_1$ and $d_2$,respectively) has the probability density
$$
f_{F(d_1,d_2)}(x) = \frac{1}{B(d_1/2,d_2/2)} (\frac{d_1}{d_2})^{d_{1}/2}x^{\frac{d_1}{2}-1}(1+\frac{d_1}{d_2}x)^{-\frac{d_1+d_2}{2}}, x>0;
$$
where $B(x,y)$ is the beta function. Thus it is a Fisher z-distribution with  appropriately changed  variable $x$ and  parameters.
\end{rem}
}

\medskip
\medskip
\medskip
\textbf{References.}

\medskip
Ph. Biane and M. Yor (1987), \emph{Variations sur une formula de Paul L\'evy},
Ann. Inst. H. Poncare Probab. Statist. 23(2), pp. 359-377.

\medskip
R. C. Bradley and Z. J. Jurek(2014), \emph{The strong mixing and the selfdecomosability},
Stat. Probab. Letters, 84, pp. 67-71.

\medskip
A. Czyzewska-Janowska and Z. J. Jurek(2011), Factorization property of generalized s-selfdecomposable
measures and the class $L^f$ distributions, \emph{Theor. Probab. Appl.} \textbf{55}, no 4, 2011, pp. 692-698.

\medskip
B. Duplantier (1989), \emph{Areas of planar Brownian curves,} J. Phys. A: Math. Gen. 22, pp. 3033-3048.

\medskip
W. Feller (1966),\emph{An introduction to probability theory and its applications}, John Wiley $\&$ Sons.

\medskip
R. K. Getoor (1979), \emph{The Brownian  escape process}, Ann. Probab. \textbf{7}, pp. 864-867.

\medskip
B.V. Gnedenko and A.N. Kolmogorov (1954), \emph{Limit distributions for sums of independent random variables}, Addison-Wesley, Reading,
Massachusetts.

\medskip
I. S. Gradshteyn and I. M. Ryzhik (1994), \emph{Tables of integrals, series, and products}, Academic Press, San Diego, Fifth Ed.

\medskip
C. Halgreen (1979), \emph{Self decomposability of the generalized inverse Gaussian and hyperbolic distributions},
 Z. Wahrsch. verw. Gebiete, 47, pp. 13-17.

\medskip
A. M. Iksanov, Z. J.Jurek, B. M. Schreiber (2004), \emph{A new factorization property of the selfdecomposable probability measures},
Ann Probab. 32, No2, 1356-1369.

\medskip
Z. J. Jurek (1985), \emph{Relations between the s-selfdecomposable ans selfdecomposable measures}, Ann. Probab. 13, 592-608.

\medskip
Z. J, Jurek (1996), \emph{Series of independent exponential random variables}, Proc.$7^{th}$ Japan-Russia Symp. Probab. Theor. Math. Stat.,
Tokyo,26-30 July 1995, World Scientific, Singapore,London, Hong Kong, 1996, pp.174- 182.

\medskip
Z. J. Jurek (1997), \emph{ Selfdecomposability: an exception or a rule?}, Annales Univer. M. Curie-Sklodowska, Sectio A, Mathematica, vol.\textbf{51}, pp. 93-107.

\medskip
Z.J. Jurek (2001), \emph{Remarks on the selfdecomposability and new examples}, Demonstratio Math. XXXIV, pp. 241-250.

\medskip
Z. J. Jurek (2003), \emph{Generalized L\'evy stochastic areas and selfdecomposability}, Stat.$\&$ Probab. Letters 64,pp.213-222.

\medskip
Z. J. Jurek (2021), \emph{Background driving distribution functions and series representations for log-gamma selfdecomposable random variables},
Theor. Probab. Applications; accepted.  also in math.arXiv: 1904.04160

\medskip
Z. J. Jurek and K. Kepczynski (2021), \emph{Graphs of the background driving distributions (BDDF) for some selfdecomposable variables,} (work in progres)

\medskip
Z. J. Jurek and J. D. Mason (1993), \emph{Opertor-limit distributions in the probability theory}, J. Wiley$\&$ Sons, New York.

\medskip
Z. J. Jurek and W. Vervaat (1983), \emph{An integral representation for selfdecomposable Banach space valued random variables}, Z. Wahrsch. verw. Gebiete 62, pp. 247-262.

\medskip
P. L\'evy (1950), \emph{ Wiener's random function and other Laplacian random functions}, Proc. $2^{nd}$ Berkeley Symp. Math. Stat. Probab., vol. II,
Berkeley CA: University of California Press, pp. 171-186.

\medskip
M. Loeve (1963), \emph{Probability theory},Third Edition, D. Van Nostrand Co., Princeton, New Jersey.

\medskip
K. Urbanik (1992), \emph{Functionals on transient stochastic processes with independent increments},  Studia Math. 103 (3), pp. 299-315.

\medskip
 N. G. Ushakov  (1999), \emph{Selected Topics in Characteristic Functions},  VSP, Utrecht, The Netherlands.

\medskip
M. Wenocur (1986),\emph{Brownian motion with quadratic killing and some implications,} J. App. Proba. 23, 893-303.

\medskip
M. Yor (1989). \emph{ On Stochastic areas and averages of plana Brownaian motion,} J. Phys. A: math. Gen 22, pp. 3049-3057

\medskip
M. Yor (1992), \emph{Some aspects of Brownian Motion}, Part I: Some special functionals, Birkhauser Verlag .

\medskip
M. Yor (1997), \emph{Some aspects of Brownian Motion}, Part II: Some special functionals, Birkhauser Verlag .

\end{document}